%% file: author_new.tex
\documentclass[graybox]{svmult}

\usepackage{mathptmx}       
\usepackage{helvet}         
\usepackage{courier}        
\usepackage{type1cm}        
\usepackage{makeidx}         
\usepackage{graphicx}        
\usepackage{multicol}        
\usepackage{multirow}        
\usepackage[bottom]{footmisc}
\usepackage{tikz} 
\usepackage{wrapfig}
\usepackage{caption} 
\usepackage{subcaption} 
\usepackage{gensymb} 
\usepackage{amsmath} 
\usepackage{amssymb} 
\usepackage{bm} 
\usepackage{comment} 
\usepackage{xcolor} 

\makeindex             

\begin{document}

\title*{Thermomechanical modelling for industrial applications }
\titlerunning{FEM based axisymmetric thermomechanical model}
\author{Nirav Vasant Shah, Michele Girfoglio and Gianluigi Rozza}
\institute{Nirav Vasant Shah, Michele Girfoglio, Gianluigi Rozza \at Scuola Internazionale Superiore di Studi Avanzati, via Bonomea, 265 - 34136 Trieste ITALY, \email{shah.nirav@sissa.it,michele.girfoglio@sissa.it,gianluigi.rozza@sissa.it}
}
%
%
\maketitle

\abstract*{In this work we briefly present a thermomechanical model that could serve as starting point for industrial applications. We address the non-linearity due to temperature dependence of material properties and heterogeneity due to presence of different materials. 
Finally a numerical example related to the simplified geometry of blast furnace hearth walls is shown with the aim of assessing the feasibility of the modelling framework.}

\abstract{In this work we briefly present a thermomechanical model that could serve as starting point for industrial applications. We address the non-linearity due to temperature dependence of material properties and heterogeneity due to presence of different materials. 
Finally a numerical example related to the simplified geometry of blast furnace hearth walls is shown with the aim of assessing the feasibility of the modelling framework.}

\section{Introduction}\label{225sec:introduction}

Thermomechanical models are widely used in many practical applications \cite{sensor_placement_paper,friction_stir_welding_paper}. We refer to the case of one-way coupling between thermal and mechanical fields, where the temperature can be computed \emph{in advance}, it being independent of the displacement, and used \emph{afterwards} to compute the displacement. Finite Element Method (FEM) \cite{fem_book} is adopted to obtain these fields by solving the weak formulation of the governing equations. FEM based thermomechanical models have been successfully used for the investigation of thermomechanical phenomena arising in blast furnace \cite{sara_paper,our_paper,phd_thesis}. In this work, we are going to take a step forward with respect to what done in \cite{our_paper} by considering the temperature dependence of material properties (that introduces a nonlinearity in the thermal model) and presence of different materials (at which one could refer to as heterogeneous material). 

\section{Mathematical model}\label{225sec:mathematical_model}


The blast furnace hearth is made up of several zones: ceramic cup, carbon block, steel shell. Each zone has different design requirement depending on the type of environment to which it is exposed. Ceramic cup is required to withstand high temperature due to direct contact with the molten metal. Carbon blocks are expected to reduce accumulation of excess heat. Steel shell is required to have sufficient mechanical strength to sustain the forces from other components. The reader is referred to \cite{our_paper,phd_thesis} for an illustration of the general layout of a blast furnace.

At the aim to consider a structure constructed using assembly of different materials, we refer to a domain $\omega$ divided into different $n_{su}$ non-overlapping subdomains $\lbrace \omega_i \rbrace_{i=1}^{n_{su}}$: 
\begin{equation}\label{225_subdomain_division}
\bar{\omega} = \bigcup\limits_{i=1}^{n_{su}} \bar{\omega}_i \ , \ \omega_i \cap \omega_j = \emptyset \ , \ i \neq j \ .
\end{equation}


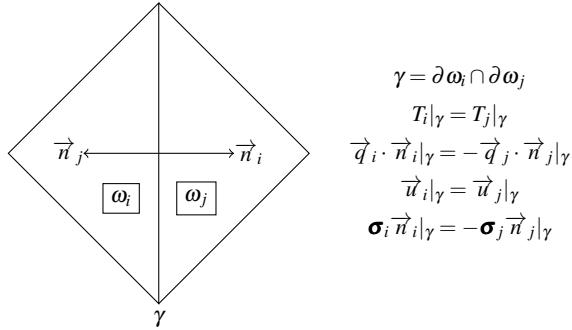
\begin{figure}
\centering
\begin{tikzpicture}
\draw (0,-2) -- (-2,0) -- (0,2) -- cycle;
\draw (0,2) -- (2,0) -- (0,-2);
\draw[->] (0,0) -- (1,0);
\draw[->] (0,0) -- (-1,0);
\node at (1.2,0) {$\overrightarrow{n}_i$};
\node at (-1.2,0) {$\overrightarrow{n}_j$};
\node[draw] at (-0.5,-0.6) {$\omega_i$};
\node[draw] at (0.5,-0.6) {$\omega_j$};
\node at (0,-2.2) {$\gamma$};
\node at (4,1) {$\gamma = \partial \omega_i \cap \partial \omega_j$};
\node at (4,0.5) {$T_i |_{\gamma} = T_j |_{\gamma}$};
\node at (4,0) {$\overrightarrow{q}_i \cdot \overrightarrow{n}_i |_{\gamma} = - \overrightarrow{q}_j \cdot \overrightarrow{n}_j |_{\gamma}$};
\node at (4,-0.5) {$\overrightarrow{u}_i |_{\gamma} = \overrightarrow{u}_j |_{\gamma}$};
\node at (4,-1) {$\bm{\sigma}_i \overrightarrow{n}_i |_{\gamma} = - \bm{\sigma}_j \overrightarrow{n}_j |_{\gamma} $};
\end{tikzpicture}
\caption{Close-up view of two subdomains (left) and continuity conditions through the interface between the two subdomains (right).}
\label{225_subdomain_interface_condition}
\end{figure}

We refer to the interface between two subdomains $\gamma = \partial \omega_i \cap \partial \omega_j \ , \ i \neq j$ as shown in Fig.  \ref{225_subdomain_interface_condition} (left). The subdomains $\omega_i$ and $\omega_j$ are related to different materials. The temperature $T$ and the heat flux $\overrightarrow{q} \cdot \overrightarrow{n}$, as well as the displacement $\overrightarrow{u}$ and the stress vector $\bm{\sigma} \overrightarrow{n}$ are continuous along the interface $\gamma$ as reported in Fig. \ref{225_subdomain_interface_condition} (right). 
Let $k^{(i)},E^{(i)},\nu^{(i)},\alpha^{(i)}$ respectively be the temperature dependent thermal conductivity, Young's modulus, Poisson's ratio and thermal expansion coefficient corresponding to the material of the subdomain $\omega_i$. In the current analysis, we consider $\nu^{(i)}$ and $\alpha^{(i)}$ constant with respect to the temperature. So, for $x \in \omega_i$, we have:
\begin{equation*}
k(T,x) = k^{(i)}(T) \ , \ E(T,x) = E^{(i)}(T) \ , \ \nu(T,x) = \nu^{(i)} \ , \ \alpha(T,x) = \alpha^{(i)} \ .
\end{equation*}
We use the piecewise spline interpolation \cite{numerical_analysis_book} to approximate thermal conductivity and Young's modulus based on their estimates (typically experimental data) related to certain discrete temperature values:
\begin{equation*}
\begin{split}
\text{if } T_a \leq T \leq T_b \ , \  k^{(i)}(T) = a_{0,k}^{(i)} T^2 + b_{0,k}^{(i)} T + c_{0,k}^{(i)} \ , \  E^{(i)}(T) = a_{0,E}^{(i)} T^2 + b_{0,E}^{(i)} T + c_{0,E}^{(i)} \ , \\
\text{if } T_b \leq T \leq T_c \ , \  k^{(i)}(T) = a_{1,k}^{(i)} T^2 + b_{1,k}^{(i)} T + c_{1,k}^{(i)} \ , \  E^{(i)}(T) = a_{1,E}^{(i)} T^2 + b_{1,E}^{(i)} T + c_{1,E}^{(i)} \ .
\end{split}
\end{equation*}
We consider a thermomechanical problem described in the cylindrical coordinate system $(r,y,\theta)$. In many real-world applications the variation of domain geometry as well as loads and heat fluxes with respect to the angular coordinate $\theta$ could be neglected. Under such conditions, it is reasonable to apply axisymmetric hypothesis. 
Then in the absence of source terms the energy and momentum conservation equations in strong formulation endowed with proper boundary conditions referred to a domain $\omega$ can be stated as follows:

\begin{subequations}\label{225_energy_equation_axisymmetric_heterogeneous}
\begin{flalign}
\text{Thermal model} & : - \frac{1}{r} \frac{\partial}{\partial r} \left( r k \frac{\partial T}{\partial r} \right) - \frac{\partial }{\partial y} \left( k \frac{\partial T}{\partial y} \right) = 0 \ , \ \text{in} \ \omega \ , \label{225_thermal_model} \\
\text{Neumann boundary} & : \left( -k(T,x) \nabla T \right) \cdot \overrightarrow{n} = 0 \ \text{on} \ \Gamma^T_N \subset \partial \omega \ , \label{225_neumann_bc} \\
\text{Convection boundary} & : \left( -k(T,x) \nabla T \right) \cdot \overrightarrow{n} = h (T - T_R) \ \text{on} \ \Gamma^T_R \subset \partial \omega \ . \label{225_convection_bc}
\end{flalign}
\end{subequations}

\begin{subequations}\label{225_momentum_equation_axisymmetric_heterogeneous}
\begin{flalign}
\text{Mechanical model} & : \frac{\partial \sigma_{rr}}{\partial r} + \frac{\partial \sigma_{ry}}{\partial y} + \frac{\sigma_{rr} - \sigma_{\theta \theta}}{r} = 0 \ , \ \text{in $\omega$} \ , \\
 & \ \ \frac{\partial \sigma_{ry}}{\partial r} + \frac{\partial \sigma_{yy}}{\partial y} + \frac{\sigma_{ry}}{r} = 0 \ , \ \text{in $\omega$} \ , \label{225_mechanical_model} \\
\text{Applied force} & : \bm{\sigma} \overrightarrow{n} = \overrightarrow{g} \ , \ \text{on} \ \Gamma^u_N \subset \partial \omega \ , \label{225_bc_applied_force} \\
\text{Bilateral frictionless contact} & : \overrightarrow{u} \cdot \overrightarrow{n} = 0 \ , \ \overrightarrow{\sigma}_t = \overrightarrow{0} \ ,  \  \text{on} \  \Gamma^u_c \subset \partial \omega \ . \label{225_bc_bilateral_frictionless_contact}
\end{flalign}
\end{subequations}

The temperature field $T$ and displacement field $\overrightarrow{u} = [u_r \ u_y]$ are the unknown quantities of interest. Convection coefficient $h$, convection temperature $T_R$ and boundary force $\overrightarrow{g}$ are specified data. The relevant material properties include thermal conductivity $k$, Young's modulus $E$, Poisson's ratio $\nu$ and thermal expansion coefficient $\alpha$. The normal vector $\overrightarrow{n}$ is considered to be pointing outwards. The shear stress $\overrightarrow{\sigma}_t$ is related to the stress tensor $\bm{\sigma}$ as:
\begin{equation*}
\overrightarrow{\sigma}_t = \bm{\sigma} \overrightarrow{n} - \sigma_n \overrightarrow{n} \ , \ \text{where} \ \sigma_n = (\bm{\sigma} \overrightarrow{n}) \cdot \overrightarrow{n} \ .
\end{equation*}

If $T_0$ is the known reference temperature, axisymmetric stress-strain relationship, in vector notation, can be expressed as,
\begin{equation*}
\lbrace \bm{\sigma}(\overrightarrow{u})[T] \rbrace = \bm{C} \lbrace \bm{\varepsilon}(\overrightarrow{u}) \rbrace - \frac{E}{(1-2\nu)} \alpha (T - T_0) \lbrace \bm{I} \rbrace \ ,
\end{equation*}
where $\bm{\varepsilon}$ is the strain tensor, $\bm{I}$ is the identity matrix and  
\begin{equation*}
\bm{C} = \frac{E}{(1 - 2 \nu)(1 + \nu)}
\begin{bmatrix}
1-\nu & \nu & \nu & 0 \\
\nu & 1-\nu & \nu & 0 \\
\nu & \nu & 1-\nu & 0 \\
0 & 0 & 0 & \frac{1 - 2 \nu}{2} \\
\end{bmatrix} \ .
\end{equation*}

We introduce weighted Sobolev spaces, $L^2_r(\omega)$ and $H^1_r(\omega)$ \cite{axisymmetric_polygonal}:
\begin{flalign*}
\begin{split}
L^2_r(\omega) & = \bigg\lbrace \psi : \omega \mapsto \mathbb{R} \ , \ \int_{\omega} \psi^2 r dr dy < \infty \bigg\rbrace \ , \\
H^1_r(\omega) & = \bigg\lbrace \psi : \omega \mapsto \mathbb{R} \ , \ \int_{\omega} \left(\psi^2 + \left( \frac{\partial \psi}{\partial r} \right)^2 +  \left( \frac{\partial \psi}{\partial y} \right)^2 \right) r dr dy < \infty \bigg\rbrace \ ,
\end{split}
\end{flalign*}
and the functional spaces for temperature and displacement:
\begin{flalign*}
\begin{split}
\mathbb{T} & = \lbrace \psi \in L^2_r(\omega) \cap H^1_r(\omega_i) \rbrace \ , \\
\mathbb{U} & = \lbrace \overrightarrow{\phi} \in [L^2_r(\omega)]^2 \ , \varepsilon(\overrightarrow{\phi}) \in [L^2_r(\omega_i)]^{3 \times 3} \ , \ \overrightarrow{\phi}\cdot\overrightarrow{n} = 0 \ \text{on $\Gamma^u_c$} \rbrace \ .
\end{split}
\end{flalign*}
Then the weak formulations corresponding to equations \eqref{225_energy_equation_axisymmetric_heterogeneous} and \eqref{225_momentum_equation_axisymmetric_heterogeneous} are given by:
\begin{align}\label{225_weak_thermal}
\sum\limits_{i=1}^{n_{su}} \int_{\omega_i} k \nabla T : \nabla \psi r dr dy & + \int_{\Gamma^T_R} h T \psi r dr dy = \int_{\Gamma^T_R} h T_R \psi r dr dy \ , \ \forall \psi \in \mathbb{T} \ ,
\end{align}
\begin{flalign}\label{225_weak_mechanical}
\begin{split}
\sum\limits_{i=1}^{n_{su}} \int_{\omega_i} \bm{C} \lbrace \bm{\varepsilon}(\overrightarrow{u}) \rbrace : \lbrace \bm{\varepsilon}(\overrightarrow{\phi}) \rbrace r dr dy  & = \sum\limits_{i=1}^{n_{su}} \int_{\omega_i} \bm{C} (T-T_0) \alpha \lbrace \bm{I} \rbrace : \lbrace \bm{\varepsilon}(\overrightarrow{\phi}) \rbrace r dr dy \\
& + \int_{\Gamma^u_N} \overrightarrow{\phi} \cdot \overrightarrow{g} r dr dy \ , \ \forall \overrightarrow{\phi} \in \mathbb{U} \ .
\end{split}
\end{flalign}

\section{Numerical example}\label{225sec:numerical_example}

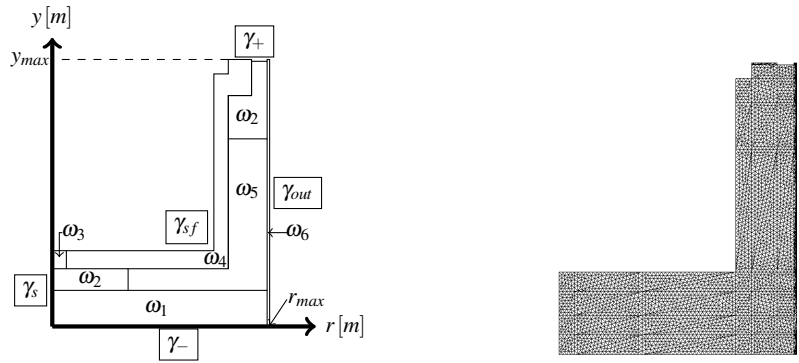
\begin{figure}
\centering
\begin{subfigure}[t]{0.52\textwidth}
\begin{tikzpicture}[scale=0.48][
dot/.style={
  fill,
  circle,
  inner sep=2pt
  }
]
\draw (5.9501,0) -- (0,0) -- (0,1.0) -- (5.9501,1.0) -- (5.9501,0);;
\node at (2.9,0.5) {$\omega_1$};
\draw (0,1.0) -- (0,1.6) -- (2.1,1.6) -- (2.1,1.0);
\node at (1.05,1.3) {$\omega_2$};
\draw (2.1,1.6) -- (4.875,1.6) -- (4.875,5.2) -- (5.9501,5.2) -- (5.9501,1);
\draw (4.875,5.2) -- (4.875,6.4) -- (5.5188,6.4) -- (5.5188,7.35) -- (5.9501,7.35) -- (5.9501,5.2);
\node at (5.4,5.8) {$\omega_2$};
\draw (0.39,1.6) -- (0.39,2.1) -- (0,2.1) -- (0,1.6);
\node (A) at (0.195,1.7) {};
\node (B) at (0.5,2.6) {};
\node at (0.6,2.6) {$\omega_3$};
\draw[->, to path={-| (\tikztotarget)}]
  (B) edge (A) ;
\draw (0.39,2.1) -- (4.475,2.1) -- (4.475,7.0) -- (4.875,7.0) -- (4.875,7.4) -- (5.5188,7.4) -- (5.5188,6.4) -- (4.875,6.4) -- (4.875,1.6);
\node at (4.52,1.85) {$\omega_4$};
\draw (5.9501,7.35) -- (5.9501,7.4) -- (6.0201,7.4) -- (6.0201,0);
\draw[->] (6.5,2.6) -- (5.97,2.6);
\node at (6.8,2.6) {$\omega_6$};
\node[draw,align=right] at (6.8,3.7) {$\gamma_{out}$};
\node[draw,align=right] at (3.5,-0.5) {$\gamma_{-}$};
\node[draw,align=right] at (3.725,2.75) {$\gamma_{sf}$};
\node[draw,align=right] at (5.6,7.9) {$\gamma_{+}$};
\node[draw,align=right] at (-0.6,1) {$\gamma_{s}$};
\node at (5.4,3.75) {$\omega_5$};
\draw[->,ultra thick] (0,0)--(7.3,0) node[right]{$r \left[m \right]$};
\draw[->,ultra thick] (0,0)--(0,8.0) node[above]{$y \left[m \right]$};
\draw[dashed] (4.875,7.4) -- (0,7.4);
\draw[->] (6.5,0.75) -- (6.0201,0.03);
\node at (7.05,0.75) {$r_{max}$};
\node at (-0.6,7.4) {$y_{max}$};
\end{tikzpicture}
\end{subfigure}
    \hfill
\begin{subfigure}[t]{0.46\textwidth}
    \centering
    \includegraphics[height=1.6in]{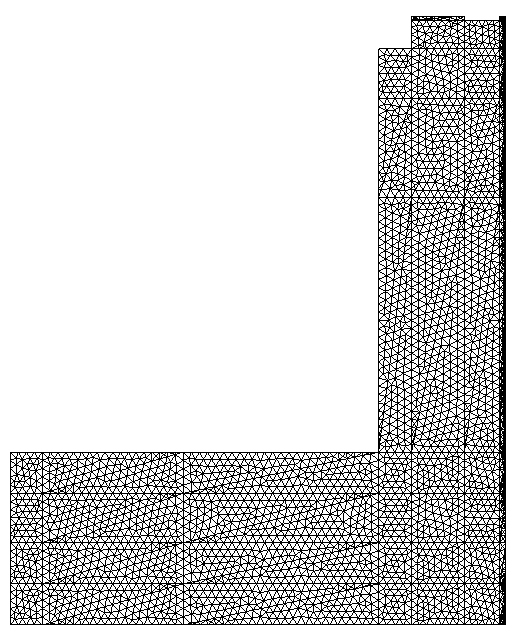}
\end{subfigure}
\caption{Computational domain (left) and view of the mesh (right).}
\label{225_hearth_2d_domain_heterogeneous}
\end{figure}

We consider the domain $\omega$ as shown in Fig. \ref{225_hearth_2d_domain_heterogeneous}. It is divided in $n_{su}=6$ subdomains. The coordinates of their vertices are reported in Table \ref{225_subdomain_coordinates}. 

At top boundary $\gamma_+$ and symmetry boundary $\gamma_s = \partial \omega \cap (r = 0)$, Neumann boundary \eqref{225_neumann_bc} and bilateral frictionless contact \eqref{225_bc_bilateral_frictionless_contact} are applied. At bottom boundary $\gamma_-$, convection boundary \eqref{225_convection_bc} and bilateral frictionless contact \eqref{225_bc_bilateral_frictionless_contact} are applied. Inner boundary $\gamma_{sf}$ and outer boundary $\gamma_{out}$ are convection and applied force boundaries (equations \eqref{225_convection_bc},\eqref{225_bc_applied_force}). Convection coefficient $h$, convection temperature $T_R$ and applied force $\overrightarrow{g}$ are reported in Table \ref{225_bc_parameters}. 

From a physical viewpoint, Neumann boundary \eqref{225_neumann_bc} on $\gamma_+$ refers to the adiabatic condition. On $\gamma_{sf}$, the convection boundary \eqref{225_convection_bc} refers to heat transfer with liquid iron at melting point. On $\gamma_{out}$ and $\gamma_-$, the convection boundary \eqref{225_convection_bc} refers to heat extraction from the structure using heat exchanger. The convection coefficients $h$ and the convection temperatures $T_R$, referring to the heat exchanger operating conditions, are kept constant. Bilateral frictionless contact \eqref{225_bc_bilateral_frictionless_contact} on $\gamma_+$ and $\gamma_-$ is related to no shear force from other components and restriction on normal expansion. The restriction on normal expansion on $\gamma_+$ refers to direct contact with other sections of hearth, while the restriction on normal expansion on $\gamma_-$ refers to direct contact with the ground. On the inner boundary $\gamma_{sf}$, the applied forces (equation \eqref{225_bc_applied_force}) refer to hydrostatic force from molten iron. Considering that the maximum hydrostatic force is exerted when the level of molten iron is $y_{max}$, we take into account the worst case scenario. 
On boundary $\gamma_{out}$, no known force occurs. 

Table \ref{225_interpolated_values} reports thermal conductivity and Young's modulus values used for the interpolation. It should be noted that the values reported in Table \ref{225_interpolated_values} are typical for the blast furnace hearth materials. The exact values of material properties depend on the commercial grade of the material used in the final design. Table \ref{225_k_E_values}\footnote{The coefficients in Table \ref{225_k_E_values} are rounded-off to maximum two decimal points.} shows the interpolation coefficients for thermal conductivity and Young's modulus. 
Thermal conductivity $k^{(6)}$ and Young's modulus $E^{(6)}$ are related to thin section of steel shell where the temperature variation is not significant, so we consider them constant. On the other hand, thermal conductivities $k^{(3)}$ and $k^{(4)}$ refer to refractory blocks, which are in direct contact with high temperature molten metal and are required to sustain high temperature. They show little variation with respect to the temperature and hence, it is reasonable to assume them constant. Table \ref{225_nu_alpha_values} reports the Poisson's ratio and thermal expansion coefficient values. The reference temperature $T_0$ is considered as $ 300 \left[ K \right]$. 

We use Lagrange finite element with polynomial of degree $1$ for displacement and temperature. The number of degrees of freedom for temperature was $4428$ and for displacement was $8856$. We use Newton's method to solve the nonlinear thermal model \eqref{225_weak_thermal} with required residual tolerance of $1e-4$. For the mechanical model \eqref{225_weak_mechanical}, we use the lower-upper (LU) decomposition. All the simulations were performed by using FEniCS \cite{fenics}.

\begin{table}
\centering
\begin{tabular}{|c|c|c|c|c|c|c|c|c|c|c|c|c|c|c|}
\hline
\multirow{2}{*} {$\omega_1$} & $r$ & 0 & 5.9501 & 5.9501 & 2.1 & 0 & / & / & / & / & / & / & / & / \\
& $y$ & 0 & 0 & 1 & 1 & 1 & / & / & / & / & / & / & / & / \\
\hline
\multirow{2}{*} {$\omega_2$} & $r$ & 0 & 2.1 & 2.1 & 0.39 & 0 & 4.875 & 5.9501 & 5.9501 & 5.5188 & 5.5188 & 4.875 & / & / \\
& $y$ & 1 & 1 & 1.6 & 1.6 & 1.6 & 5.2 & 5.2 & 7.35 & 7.35 & 7 & 7 & / & / \\
\hline
\multirow{2}{*} {$\omega_3$} & $r$ & 0 & 0.39 & 0.39 & 0 & / & / & / & / & / & / & / & / & / \\
& $y$ & 1 & 1 & 1.6 & 1.6 & / & / & / & / & / & / & / & / & / \\
\hline
\multirow{2}{*} {$\omega_4$} & $r$ & 0.39 & 2.1 & 4.875 & 4.875 & 4.875 & 5.5188 & 5.5188 & 5.5188 & 4.875 & 4.875 & 4.475 & 4.475 & 0.39 \\
& $y$ & 1.6 & 1.6 & 1.6 & 5.2 & 6.4 & 6.4 & 7.35 & 7.4 & 7.4 & 7 & 7 & 2.1 & 2.1 \\
\hline
\multirow{2}{*} {$\omega_5$} & $r$ & 2.1 & 5.9501 & 5.9501 & 4.875 & 4.875 & 2.1 & / & / & / & / & / & / & / \\
& $y$ & 1 & 1 & 5.2 & 5.2 & 1.6 & 1.6 & / & / & / & / & / & / & / \\
\hline
\multirow{2}{*} {$\omega_6$} & $r$ & 5.9501 & 5.9501 & 5.9501 & 5.9501 & 5.9501 & 6.0201 & 6.0201 & / & / & / & / & / & / \\
& $y$ & 0 & 1 & 5.2 & 7.35 & 7.4 & 7.4 & 0 & / & / & / & / & / & / \\
\hline
\end{tabular}
\caption{Coordinates (in $\left[ m \right]$) of the vertices of subdomains $\lbrace \omega_i \rbrace_{i=1}^6$ (see Fig. \ref{225_hearth_2d_domain_heterogeneous}).}
\label{225_subdomain_coordinates}
\end{table}

\begin{table}
\centering
\begin{tabular}{|c|c|c|c|c|c|}
\hline
Boundary & $\gamma_-$ & $\gamma_{out}$ & $\gamma_+$ & $\gamma_{sf}$ & $\gamma_s$\\ 
\hline
Convection coefficient $h \left[ \frac{W}{m^2K} \right]$ & $200$ & $200$ & / & $2000$ & / \\
\hline
Convection temperature $T_R \left[ K \right]$ & $300$ & $300$ & / & $1773$ & / \\
\hline
Boundary force $\overrightarrow{g} \left[ \frac{N}{m^2} \right] $ & / & $\overrightarrow{0}$ & / & $-77106 (y_{max}-y) \overrightarrow{n}$ & / \\
\hline
\end{tabular}
\caption{Convection coefficients and temperatures, applied forces at domain boundaries.}
\label{225_bc_parameters}
\end{table}

\begin{table}
\centering
\begin{tabular}{|c|c|c|c|c|c|c|c|c|c|c|c|c|c|}
\hline
$T \left[ K \right]$ & \multicolumn{6}{|c|}{Thermal conductivity $\left[\frac{W}{mK}\right]$} & $T \left[K \right]$ & \multicolumn{6}{|c|}{Young's modulus $\left[GPa \right]$} \\
\hline
& $\omega_1$ & $\omega_2$ & $\omega_3$ & $\omega_4$ & $\omega_5$ & $\omega_6$& & $\omega_1$ & $\omega_2$ & $\omega_3$ & $\omega_4$ & $\omega_5$ & $\omega_6$\\
\hline
293 & 16.07 & 49.35 & 5.3 & 4.75 & 23.34 & 45.6 & 293 & 10.5 & 15.4 & 58.2 & 1.85 & 14.5 & 190 \\
473 & 15.53 & 24.75 & 5.3 & 4.75 & 20.81 & 45.6 & 573 & 10.3 & 14.7 & 67.3 & 1.92 & 15.0 & 190 \\
873 & 15.97 & 27.06 & 5.3 & 4.75 & 20.99 & 45.6 & 1073 & 10.4 & 13.8 & 52.9 & 1.83 & 15.3 & 190 \\
1273 & 17.23 & 38.24 & 5.3 & 4.75 & 21.62 & 45.6 & 1273 & 10.3 & 14.4 & 51.6 & 1.85 & 13.3 & 190 \\
\hline
\end{tabular}
\caption{Temperature dependent thermal conductivity and Young's modulus values used for interpolation.}
\label{225_interpolated_values}
\end{table}

\begin{table}
\centering
\begin{tabular}{|c|c|c|c|c|c|c|c|c|c|c|c|c|}
\hline
 & $a_{0,k}^{(i)}$ & $b_{0,k}^{(i)}$ & $c_{0,k}^{(i)}$ & $a_{1,k}^{(i)}$ & $b_{1,k}^{(i)}$ & $c_{1,k}^{(i)}$ & $a_{0,E}^{(i)}$ & $b_{0,E}^{(i)}$ & $c_{0,E}^{(i)}$ & $a_{1,E}^{(i)}$ & $b_{1,E}^{(i)}$ & $c_{1,E}^{(i)}$\\
\hline
$\omega_1$ & 1.4E-5 & -1.3E-2 & 1.9E1 & -4.7E-6 & 1.1E-2 & 1.1E1 & 1.2E3 & -1.6E6 & 1.1E10 & -1.2E3 & 2.4E6 & 9.2E9\\ 
\hline
$\omega_2$ & 3.9E-4 & -4.3E-1 & 1.4E2 & -1.2E-4 & 2.5E-1 & -8.7E1 & -4.5E2 & -2.3E6 & 1.6E10 & 9.1E3 & -1.8E7 & 2.3E10\\
\hline
$\omega_3$ & / & / & 5.3 & / & / & 5.3 & -1.05E5 & 1.2E8 & 3.1E10 & 6.1E4 & -1.5E8 & 1.4E11\\
\hline
$\omega_4$ & / & / & 4.75 & / & / & 4.75 & -7.4E2 & 8.8E5 & 1.7E9 & 6.5E2 & -1.4E6 & 2.6E9\\
\hline
$\omega_5$ & 3.9E-5 & -4.4E-2 & 3.3E1 & -1.3E-5 & 2.6E-2 & 9.2 & 1.3E3 & 8.9E5 & 1.4E10 & -1.9E4 & 3.4E7 & 5.6E8\\
\hline
$\omega_6$ & / & / & 45.6 & / & / & 45.6 & / & / & 1.9E11 & / & / & 1.9E11\\
\hline
\end{tabular}
\caption{Interpolation coefficients for thermal conductivity ($T_a = 293 \left[K \right], T_b = 673 \left[K \right], T_c = 1800 \left[K \right]$) and for Young's modulus ($T_a = 293 \left[K \right], T_b = 823 \left[K \right], T_c = 1800 \left[K \right]$).}
\label{225_k_E_values}
\end{table}

\begin{table}
\centering
\begin{tabular}{|c|c|c|c|c|c|c|}
\hline
 & $\omega_1$ & $\omega_2$ & $\omega_3$ & $\omega_4$ & $\omega_5$ & $\omega_6$\\
\hline
$\nu^{(i)}$ & 0.3 & 0.2 & 0.1 & 0.1 & 0.2 & 0.3 \\
\hline
$\alpha^{(i)} \left[ K^{-1} \right]$  & 2.3E-6 & 4.6E-6 & 4.7E-6 & 4.6E-6 & 6E-6 & 1.2E-5 \\
\hline
\end{tabular}
\caption{Poisson's ratio and thermal expansion coefficient values.}
\label{225_nu_alpha_values}
\end{table}

\begin{figure}
\begin{subfigure}[t]{0.48\textwidth}
    \centering
    \includegraphics[height=1.6in]{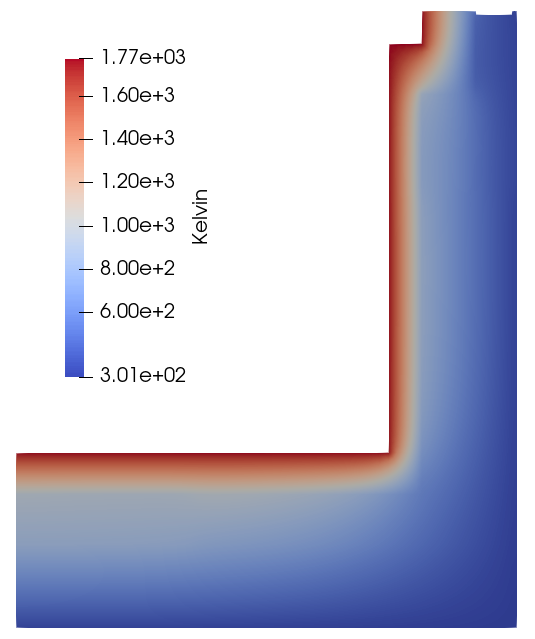}
\end{subfigure}
    \hfill
\begin{subfigure}[t]{0.48\textwidth}
    \centering
    \includegraphics[height=1.6in]{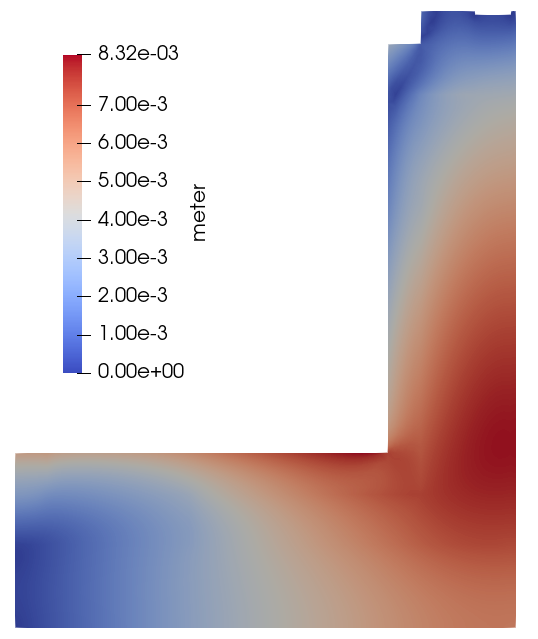}
\end{subfigure}
\caption{Computed temperature field (left) and displacement field (right).}
\label{225_computed_results}
\end{figure}

Numerical results are shown in Fig. \ref{225_computed_results}. As can be noticed, both temperature and displacement profiles do not show strong discontinuity at the interfaces. This demonstrates that the interface conditions related to heat flux and stresses (see Fig. \ref{225_subdomain_interface_condition}) are properly formulated and incorporated in the weak formulation.

From practical viewpoint, the temperature profile is typically used to identify areas subjected to high thermal stress. In addition, the temperature profile is also used to locate critical isotherms in the domain, such as isotherm corresponding to $1150 \degree C$, which represents the penetration of liquid iron in the blast furnace hearth. On the other hand, the displacement profile is typically used to identify areas with maximum deformation. Furthermore, the displacement field along with temperature field have direct impact on the stress field in the hearth.

\section{Concluding remarks}\label{225sec:concluding_remarks}

In this work we have addressed the development of a thermomechanical model able to describe phenomena associated to the temperature dependence of material properties (non linearity) and to the presence of different materials (heterogeneity). 
We expect this preliminary work could serve as starting point for thermomechanical analysis of practical problems. 

\begin{acknowledgement}
We acknowledge the financial support of the European Union under the Marie Sklodowska-Curie Grant Agreement No. 765374 and the partial support by the European Union Funding for Research and Innovation - Horizon 2020 Program - in the framework of European Research Council Executive Agency: Consolidator Grant H2020 ERC CoG 2015 AROMA-CFD project 681447 ``Advanced Reduced Order Methods with Applications in Computational Fluid Dynamics''. This work has focused exclusively on civil applications. It is not to be used for any illegal, deceptive, misleading or unethical purpose or in any military applications causing death, personal injury or severe physical or environmental damage.
\end{acknowledgement}

\input{referenc}

\end{document}

%% file: referenc.tex
%
%
%